\documentclass[12pt]{article}
\usepackage{amsfonts}
\usepackage{amsmath}
\usepackage{amssymb}
\newcommand{\proof}{\noindent{\em Proof: }}
\newcommand{\qed}{\hspace{\fill}$\square$}
\newtheorem{theorem}{Theorem}

\newtheorem{definition}{Definition}
\newtheorem{example}{Example}
\newtheorem{lemma}{Lemma}

\numberwithin{equation}{section}
\date{ }

\title{The Chaos Game on a General Iterated Function System}
\author{Michael F. Barnsley \\
Department of Mathematics \\ Australian National University \\
Canberra, ACT, Australia \\
{\tt michael.barnsley@maths.anu.edu.au} \\
http://www.superfractals.com \\ \\
Andrew Vince \\
Department of Mathematics \\ University of Florida\\
Gainesville, FL 32611-8105, USA \\
{\tt avince@ufl.edu} }

\begin{document}
\maketitle
\vskip 1.5cm 

\begin{abstract}
The main theorem of this paper establishes conditions under which the "chaos
game" algorithm almost surely yields the attractor of an iterated function
system. The theorem holds in a very general setting, even for non contractive
iterated function systems, and under weaker conditions on the random orbit of
the chaos game than obtained previously. 
\end{abstract}
\newpage
\baselineskip 20pt

\section{Introduction}
There are two methods for computing pictures of fractals that are
attractors of iterated functions systems, the deterministic algorithm
and the more efficient ``chaos game'' algorithm \cite{barnsley}.  This
paper concerns the chaos game on a general iterated function system
(IFS) $\mathcal{F}$ defined on a complete metric space
$(\mathbb{X},d)$.  The iterated function system is ``general'' in the
following sense. The only restriction placed on $\mathbb{X}$ is that
it is proper, i.e., closed balls are compact. The only restriction on
the functions in $F$ is that they be continuous. In particular, they
need not be contractions on $\mathbb{X}$. In fact, none of the
functions needs to be a contraction on $\mathbb{X}$ with respect to
any metric giving the same topology as the original metric $d$.  This
paper provides a natural definition of an attractor of such a general
IFS.  A general IFS may possess more than one attractor, one
attractor, or no attractor. Examples of iterated function systems that are non
contractive yet possess attractors are given in
Section~\ref{examples}.

The main result, Theorem \ref{chaosthm}, is new in that it shows that the
chaos game algorithm can be applied to such a general IFS - in particular to
situations where there is no metric, even in a neighborhood of the attractor,
with respect to which the IFS is contractive. We show that the chaos game
algorithm almost always yields the attractor. More precisely, if the IFS has
an attractor $A$, then a random orbit, starting with a point in the basin of
attraction of $A$, converges with probability one to the attractor. We show
this under weaker conditions than have heretofore been described. In all other
papers on this topic, of which we are aware, for example \cite{barnsleydemko},
\cite{barnplacedep}, \cite{berger}, \cite{JEergodic}, \cite{jaroszewska},
\cite{vrscay1}, \cite{onicescu}, \cite{werner}, \cite{S}, it is required that
the IFS be either contractive, or contractive on the average. It is also
required that the process by which the functions are selected to generate the
orbit is stationary and that the selection process depends H\"{o}lder
continuously on the initial point; see for example \cite{stenflo3}. For our
result, none of these conditions are required. The only constraint on how the
functions are randomly selected, one after another, is that, at each step,
there is the possibility of choosing each and any of the functions, and that
the probability of so doing is bounded below away from zero.

We do not ask, as is customary in descriptions of the chaos game algorithm,
that the process be i.i.d. or that it be conditioned on the "past" in some
restricted way. This is because we are not concerned with the existence of an
invariant measure associated with a chaos game on an IFS. For such an
invariant measure to exist it is typically required that an associated
transition probability operator is suitably well-behaved, see \cite{jamison}
and references therein. Since we are only \ concerned with the relationship
between random orbits and attractors of an IFS, we are able to obtain almost
sure convergence of random orbits to attractors in very general situations.

Section 2 of this paper contains basic definitions, in particular the
definition of an attractor of an IFS. Theorem~\ref{attractorthm} in Section 2
provides an expression for an attractor, of some independent interest, that
will be used to prove the main result. The main result on the chaos game is
Theorem \ref{chaosthm} in Section 3. Section 4 contains examples that
illustrate the practical value of Theorem \ref{chaosthm}.

\section{\label{general}General iterated function systems}

Throughout this paper $(\mathbb{X},d_{\mathbb{X}})$ is a complete metric space.

\begin{definition} \label{defIFS}
If $f_{m}:\mathbb{X}\rightarrow\mathbb{X}$, $m=1,2,\dots,M,$ are continuous
mappings, then $\mathcal{F}=\left(  \mathbb{X};f_{1},f_{2},...,f_{M}\right)  $
is called an \textbf{iterated function system} (IFS).
\end{definition}

By slight abuse of terminology we use the same symbol $\mathcal{F}$ for the
IFS, the set of functions in the IFS, and for the following mapping. Letting
$2^{\mathbb{X}}$ denote the collection of subsets of $\mathbb{X}$, define
$\mathcal{F}:2^{\mathbb{X}}\mathbb{\rightarrow}2^{\mathbb{X}}$ by%
\[
\mathcal{F}(B)=\bigcup_{f\in\mathcal{F}}f(B)
\]
for all $B\in2^{\mathbb{X}}$.

Let $\mathbb{H=H(X)}$ be the set of nonempty compact subsets of $\mathbb{X}$.
Since $\mathcal{F}\left(  \mathbb{H}\right)  \subset\mathbb{H}$ we can also
treat $\mathcal{F}$ as a mapping $\mathcal{F}:\mathbb{H\rightarrow H}$. Let
$d_{\mathbb{H}}$ denote the Hausdorff metric on $\mathbb{H}$, defined in terms
of $d := d_{\mathbb{X}}$. Using the notation
\[
S+r = \{y\in\mathbb{X} \, : \, d(x,y)<r \text{ for some} x\in S\}
\]
with $S\subset\mathbb{X}$ and $r>0,$ a convenient definition of the Hausdorff
metric (see for example \cite[p.66]{edgar}) is
\[
d_{\mathbb{H}}(B,C)=\inf\{r>0:B\subset C+r \text{ and } C\subset B+r\}
\]
for all $B,C\in\mathbb{H}$.

A metric space $\mathbb{X}$ is \textit{locally compact} if every point has a
compact neighborhood and is \textit{proper} if every closed ball $\{y \, : \,
d(x,y) \leq r\}$ is compact. Proper spaces are locally compact, but the
converse is not true in general.

\begin{lemma}
\label{proper}

\begin{enumerate}
\item A metric space is proper if and only if $C+r$ is compact whenever
$C\subset\mathbb{X}$ is compact and $r$ is a positive real number.

\item If $\mathbb{X}$ is proper, then $\mathcal{F}:\mathbb{H }(\mathbb{X})
\rightarrow\mathbb{H }(\mathbb{X})$ is continuous.
\end{enumerate}
\end{lemma}

\proof
Concerning statement (1), clearly the condition on $C+r$ in the statement of
the lemma implies that the metric space is proper. Just take the set $C$ to be
a single point.

Conversely assume that the metric space $\mathbb{X}$ is proper. Because $C$ is
compact, $C$ is totally bounded. This implies that $C$ is bounded, and hence
$C+r$ is bounded. Therefore there is a closed ball $B$ such that $C+r \subset
B$. But $B$ is compact because $\mathbb{X}$ is proper. Therefore $B$ is
totally bounded. This implies that, for any $\varepsilon>0$, there is a finite
set $\{ B_{i} := B(x_{i},\varepsilon)\, : \, i = 1,2, \dots, q \}$ of closed
balls of radius $\varepsilon$ centered at the points $x_{i}$ such that $C+r
\subset B \subset\cup_{i=1}^{q} B_{i}$. Therefore $C+r$ is totally bounded.
Since $C+r$ is also closed, it is compact.

Concerning statement (2), for any $B\in\mathbb{H}$, we will show that
$\mathcal{F}:\mathbb{H\rightarrow H}$ is continuous at $B$. Since $\mathbb{X}$
is proper and $B$ is compact, statement (1) implies that $B+1$ is compact.
Therefore each $f\in\mathcal{F}$ is uniformly continuous on $B+1$. It follows
that, for any $\varepsilon> 0$, there is a $\delta>0$ such that $d_{\mathbb{X}%
}(f(x),f(y))<\varepsilon$ whenever $d_{\mathbb{X}}(x,y)<\delta,$ for all
$x,y\in$ $B+1$ and for all $f\in\mathcal{F}.$ Without loss of generality, take
$\delta<1$.

Let $C\in\mathbb{H}$ with $d_{\mathbb{H}}(B,C)<\delta$ and let $f\in
\mathcal{F}$. We will show that $d_{\mathbb{H}}(f(B),f(C))<\varepsilon$. Let
$b^{\prime}\in f(B)$ and $b\in B$ such that $f(b)=b^{\prime}$. Since
$d_{\mathbb{H}}(B,C)<\delta$ there is $c\in C$ such that $d(b,c)<\delta$.
Since $\delta_{\varepsilon}<1$ we have $c\in B+1$. It follows that
$d(f(b),f(c))<\varepsilon$ and therefore that $f(B)\subset f(C)+\varepsilon.$
By a similar argument $f(C)\subset f(B)+\varepsilon.$ Hence $d_{\mathbb{H}
}(f(B),f(C))<\varepsilon$ for all $f \in\mathcal{F}$. Finally $d_{\mathbb{H}}(
\mathcal{F} (B),\mathcal{F} (C)) \leq\max_{f\in\mathcal{F}} d_{\mathbb{H}} (
f(B), f(C) )$.
\qed

Statement (i) of the following foundational results is well-known. 
A short proof
can be found in \cite[p.67, Theorem 2.4.4]{edgar}. Statement (ii), also well
known, can be found, for example, in \cite{henrikson}. Statement (iii) is a
classical result of Hutchinson \cite{hutchinson}.

\begin{theorem}
\label{ctythm}(i) The metric space $(\mathbb{H},d_{\mathbb{H}})$ is complete.

(ii) If $(\mathbb{X},d_{\mathbb{X}})$ is compact then $(\mathbb{H}%
,d_{\mathbb{H}})$ is compact.

(iii) If $f:\mathbb{X\rightarrow}\mathbb{X}$ is a contraction mapping for each
$f\in\mathcal{F}$, then $\mathcal{F}:\mathbb{H\rightarrow H}$ is a 
contraction mapping.
\end{theorem}

For $B\subset\mathbb{X}$, let $\mathcal{F}^{k}(B)$ denote the $k$-fold
composition of $\mathcal{F}$, the union of $f_{i_{1}}\circ f_{i_{2}}%
\circ\cdots\circ f_{i_{k}}(B)$ over all finite words $i_{1}i_{2}\cdots i_{k}$
of length $k.$ Define $\mathcal{F}^{0}(B)=B.$

\begin{definition}
\label{attractdef}A nonempty compact set $A\subset\mathbb{X}$ is said to be an
\textbf{attractor} of the IFS $\mathcal{F}$ if

(i) $\mathcal{F}(A)=A$ and

(ii) there is an open set $U\subset\mathbb{X}$ such that $A\subset U$ and
$\lim_{k\rightarrow\infty}\mathcal{F}^{k}(B)=A,$ for all $B\in\mathbb{H(}U)$,
where the limit is with respect to the Hausdorff metric.

The largest open set $U$ such that (ii) is true is called the \textbf{basin of
attraction} for the attractor $A$ of the IFS $\mathcal{F}$.
\end{definition}

The following observation \cite[Proposition 3 (vii)]{lesniak}, \cite[p.68,
Proposition 2.4.7]{edgar} will be used in proving Theorem~\ref{attractorthm}.

\begin{lemma}
\label{intersectlemma}Let $\left\{  B_{k}\right\}  _{k=1}^{\infty}$ be a
sequence of nonempty compact sets such that $B_{k+1}\subset B_{k}$ for all
$k$. Then $\cap_{k\geq1}B_{k}=\lim_{k\rightarrow\infty}B_{k}$ where
convergence is with respect to the Haudorff metric.
\end{lemma}

The notation $\overline{S}$ is used to denote the closure of a set $S$, and,
when $U\subset\mathbb{X}$ is nonempty, $\mathbb{H(}U)=\mathbb{H(X)\cap}2^{U}$.
The quantity on the right-hand side of the equation below is sometimes called
the \textit{topological upper limit} of the sequence $\left\{  F^{k}%
(B)\right\}  _{k=1}^{\infty}$ and is related to other definitions of
attractors of generalizations of the notion of an IFS; see for example McGehee
\cite{mcgehee} and Lesniak \cite{lesniak}, and references in both of these.
We will use Theorem
\ref{attractorthm} in the proof of Theorem \ref{chaosthm}, our main result.

\begin{theorem}
\label{attractorthm} Let $\mathcal{F}$ be an IFS with attractor $A$, with
basin of attraction $U.$ If $\mathcal{F}:\mathbb{H(}U)\mathbb{\rightarrow
H(}U)$ is continuous then
\[
A=\bigcap\limits_{K\geq1}\overline{\bigcup_{k\geq K}\mathcal{F}^{k}(B)}%
\]
for all $B\subset U$ such that $\overline{B}\in\mathbb{H(}U)$.
\end{theorem}

\proof
We carry out the proof under the assumption that $B\in\mathbb{H(}U)$. It then
follows from \cite[Proposition 3 (i)]{lesniak} that Theorem \ref{attractorthm}
is true for all $B\subset U$ such that $\overline{B}\in\mathbb{H(}U).$

Clearly $\overline{\bigcup_{k\geq K} \mathcal{F} ^{k}(B)}$ is nonempty 
because $B$ is
nonempty. We next show that $\overline{\bigcup_{k\geq K} \mathcal{F} ^{k}(B)}$ 
is compact.
Let $\left\{  \mathcal{O}_{i}\,:\,i\in\mathcal{I}\right\}  $ be an open cover
of $\overline{\bigcup_{k\geq K}\mathcal{F}^{k}(B)}$. Since $\bigcup_{k\geq K}
\mathcal{F}^{k}(\overline{B})\subset\overline{\bigcup_{k\geq K}\mathcal{F}^{k}(B)}$ and
$A=\lim_{k\rightarrow\infty}\mathcal{F}^{k}(\overline{B})$, also $\left\{
\mathcal{O}_{i}\,:\,i\in\mathcal{I}\right\}  $ is an open cover of $A$.
Because $A$ is compact, $\left\{  \mathcal{O}_{i}\,:\,i\in\mathcal{I}\right\}
$ contains a finite subcollection, say $\left\{  \mathcal{O}_{m}%
\,:\,m=1,2,..,M\right\}  $, such that $A\subset\mathcal{O}:=\cup_{m=1}
^{M}\mathcal{O}_{m}$. Because a metric space is normal, there is an open set
$\mathcal{O^{\prime}}$ containing $A$ such that $\overline{\mathcal{O^{\prime
}}}\subset\mathcal{O}$. Again using that fact that $\mathcal{F}^{k}(\overline{B})$
converges in the Hausdorff metric to $A$, there is an integer $K^{\prime}$
such that $\mathcal{F}^{k}(\overline{B})\subset\mathcal{O^{\prime}}$ for all $k\geq
K^{\prime}$. It follows that $\bigcup_{k\geq K^{\prime}}\mathcal{F}^{k}(\overline
{B})\subset\mathcal{O^{\prime}}$ and therefore $\overline{\bigcup_{k\geq
K}\mathcal{F}^{k}(B)}\subset\mathcal{O}=\cup_{m=1}^{M}\mathcal{O}_{m}$ for all $K\geq
K^{\prime}$. Therefore $\overline{\bigcup_{k\geq K}\mathcal{F}^{k}(B)}$ is compact if
$K\geq K^{\prime}$. If $K\leq K^{\prime}$, then%
\begin{equation}
\overline{\bigcup_{k\geq K}\mathcal{F}^{k}(B)}=\bigcup_{K^{\prime}>k\geq K}%
\overline{\mathcal{F}^{k}(B)}\cup\overline{\bigcup_{k\geq K^{\prime}}\mathcal{F}^{k}(B)}.
\label{Eq1}%
\end{equation}
Since each function $f\in \mathcal{F}$ is continuous and $\overline{B}$ is compact,
$\overline{\mathcal{F}^{k}(B)}=\mathcal{F}^{k}(\overline{B})$ is compact. Hence $\overline
{\bigcup_{k\geq K}\mathcal{F}^{k}(B)}$, the finite union of compact sets on the right
hand side of equation~\ref{Eq1}, is also compact. Since $\overline
{\bigcup_{k\geq K}\mathcal{F}^{k}(B)}$ is a nonempty compact set, the nested
intersection%
\[
\tilde{A}:=\bigcap\limits_{K\geq1}\overline{\bigcup\limits_{k\geq
K}\mathcal{F}^{k}(B)}%
\]
is a nonempty compact set.

Now observe that%
\begin{align*}
\mathcal{F}\left(  \tilde{A}\right)   &  =\mathcal{F}\left(  \lim
_{K\rightarrow\infty}\overline{\bigcup\limits_{k\geq K}\mathcal{F}^{k}%
(B)}\right)  \qquad\text{ (by Lemma \ref{intersectlemma})}\\
&  =\lim_{K\rightarrow\infty}\mathcal{F}\left(  \overline{\bigcup
\limits_{k\geq K}\mathcal{F}^{k}(B)}\right)  \text{ }\qquad\text{(since
}\mathcal{F}:\mathbb{H\rightarrow H}\text{ is continuous)}\\
&  =\lim_{K\rightarrow\infty}\bigcup\limits_{f\in\mathcal{F}}f\left(
\overline{\bigcup\limits_{k\geq K}\mathcal{F}^{k}(B)}\right)  \text{ }%
\qquad\text{(by the definition of }\mathcal{F}\text{)}\\
&  =\lim_{K\rightarrow\infty}\bigcup\limits_{f\in\mathcal{F}}\overline
{f\left(  \bigcup\limits_{k\geq K}\mathcal{F}^{k}(B)\right)  }\text{ (since
}f:\mathbb{X\rightarrow X}\text{ is continuous and }\\
&  \qquad\qquad\qquad\qquad\qquad\qquad\qquad\overline{\cup_{k\geq
K}\mathcal{F}^{k}(B)}\text{ is compact)}\\
&  =\lim_{K\rightarrow\infty}\overline{\bigcup\limits_{f\in\mathcal{F}%
}f\left(  \bigcup\limits_{k\geq K}\mathcal{F}^{k}(B)\right)  }\text{ }%
\qquad\text{(since }\bigcup\limits_{f\in\mathcal{F}}\text{is a finite
union)}\\
&  =\lim_{K\rightarrow\infty}\overline{\bigcup\limits_{k\geq K}\bigcup
\limits_{f\in\mathcal{F}}f(\mathcal{F}^{k}(B))}\\
&  =\lim_{K\rightarrow\infty}\overline{\bigcup\limits_{k\geq K}\mathcal{F}%
^{k}(B)}=\tilde{A}\text{.}%
\end{align*}
By an argument essentially the same as used above to show that 
\linebreak $\overline
{\bigcup_{k\geq K}\mathcal{F}^{k}(B)}\subset\mathcal{O}$, it can be shown that
$\overline{\cup_{k\geq K}\mathcal{F}^{k}(B)}$ for $K$ sufficiently large, and
hence $\tilde{A}$, lies in $U$.

Since $\tilde{A}$ is nonempty, compact and lies in $U$, the basin of
attraction of $A$, we must have
\[
A=\lim_{k\rightarrow\infty}\mathcal{F}^{k}(\tilde{A})=\tilde{A}.
\]
\qed

\section{The chaos game algorithm}

The following lemma will be used in the proof of Theorem~\ref{chaosthm}.

\begin{lemma}
\label{uniformlemma} Let $\mathbb{X}$ be proper, and let $\mathcal{F} =\left(
\mathbb{X};f_{1},f_{2},...,f_{N}\right)  $ be an IFS with attractor $A$. For
any $\varepsilon>0$ there is an integer $M=M(\varepsilon)$ such that for each
$x\in A+\varepsilon$ there is an integer $m=m(x,\varepsilon)< M$ such that
\[
d_{\mathbb{H}}(A,\mathcal{F}^{m}(\left\{  x\right\}  )<\varepsilon/2.
\]
\end{lemma}

\proof
Because $\mathbb{X}$ is proper and $A$ is compact, $A+\varepsilon$ is also
compact by statement (1) of Lemma~\ref{proper}. There is no loss of generality
in assuming that $\varepsilon$ is sufficiently small that $A+\varepsilon
\subset U$, where $U$ is the basin of attraction of $A$. If $x\in
A+\varepsilon\subset U$, then there is an integer $m(x,\varepsilon)\geq0$ such
that
\begin{equation}
\label{eq2}d_{\mathbb{H}}(A,\mathcal{F}^{m(x,\varepsilon)}(\left\{  x\right\}
)<\varepsilon/4.
\end{equation}
This is because $\lim_{k \rightarrow\infty} \mathcal{F}^{k}(\{x\}) = A$.

Since $\mathbb{X}$ is proper, it follows from statement (2) of
Lemma~\ref{proper} that $\mathcal{F}:\mathbb{H\rightarrow H}$ is continuous,
whence $\mathcal{F}^{m(x,\varepsilon)}:\mathbb{H\rightarrow H}$ is continuous.
Since $\mathcal{F}^{m(x,\varepsilon)}:\mathbb{H\rightarrow H}$ is continuous,
there is an open ball $B(\left\{  x\right\}  ,r_{x})$ (in $\mathbb{H}$) of
radius $r_{x}>0$ centered at $\{x\}$ such that $d_{\mathbb{H}}(\mathcal{F}%
^{m(x,\varepsilon)}\{x\},\mathcal{F}^{m(x,\varepsilon)}(Y))<\varepsilon/4$ for
all $Y\in B(\left\{  x\right\}  ,r_{x})$. It follows, in particular, that
there a ball $B(x,r_{x})$ (in $\mathbb{X}$) centered at $x$ such that
$d_{\mathbb{H}} (\mathcal{F}^{m(x,\varepsilon)}\{x\},\mathcal{F}
^{m(x,\varepsilon)}(\left\{  y\right\}  ))<\varepsilon/4$ for all $y\in
B(x,r_{x})$. Combining this with equation~\ref{eq2} above gives $d_{\mathbb{H}%
} (A,\mathcal{F} ^{m(x,\varepsilon)}(\left\{  y\right\}  ))<\varepsilon/2$ for
all $y\in B(x,r_{x})$.

The set of balls $\{ B(x,r_{x}), \,x \in A+\varepsilon\}$ is an open covering
of $A_{\varepsilon}$. Since $A+\varepsilon$ is compact, there is a finite
subcovering $\{ B(a_{i},r_{a_{i}})\, : \, i = 1,2, \dots, q\}$ such that
$A+\varepsilon\subset\cup_{i=1}^{q}B(a_{i},r_{a_{i}})$. If $M=\max_{i}%
m(a_{i},\varepsilon)$, then for any $x\in A+\varepsilon$ there is a $i$ such
that that $x\in B(a_{i},r_{a_{i}})$, in which case $d_{\mathbb{H}}
(A,\mathcal{F}^{m}(\left\{  x\right\}  ))<\varepsilon/2$ for
$m=m(x,\varepsilon):=m(a_{i},\varepsilon)< M$.
\qed

\begin{definition}
Let $\mathcal{F}=\left(  \mathbb{X};f _{1},f_{2},...,f_{N}\right)  $ be an IFS
and $p\in(0,1/N]$ fixed. A sequence $\left\{  x_{k}\right\}  _{k=0}^{\infty}$
of points in $\mathbb{X}$ is called a \textbf{random orbit} of $x_{0}%
\in\mathbb{X}$ if $x_{k}=f_{\sigma_{k}}(x_{k-1}), \, k\geq1$ where $\sigma
_{k}$ is selected randomly from $\{1,2,...,N\}$, for $k=1,2,...$, where the
probability that $\sigma_{k}=n$ is greater than\ or equal to $p,$ regardless
of the preceeding outcomes, for all $n\in\{1,2,...,N\}$ and all $k.$ More
formally, in terms of conditional probability,
\[
P(\sigma_{k}=n \, | \, x_{0},\sigma_{1},\sigma_{2},...,\sigma_{k-1})\geq p.
\]
\end{definition}

\begin{theorem}
\label{chaosthm} Let $\mathbb{X}$ be a proper complete metric space and
\linebreak
$\mathcal{F}=\left(  \mathbb{X};f_{1},f_{2},...,f_{N} \right)  $ an IFS with
attractor $A$ and basin of attraction $U$. If $\left\{  x_{k}\right\}
_{k=0}^{\infty}$ is a random orbit of $x_{0}\in U$ under $\mathcal{F}$, then
with probability one,
\[
A=\lim_{K\rightarrow\infty} \left\{ x_{k}\right\}_{k=K}^{\infty},\]
where the limit is with respect to the Hausforff metric.
\end{theorem}

\proof We first claim that, given any $\varepsilon>0$, there is an 
integer $K>0$ such that 
\begin{equation} \label{eq} x_{k} \in A+\varepsilon \end{equation}
 for all $k \geq K$.
Since $\mathbb{X}$ is proper, $\mathcal{F}$ is continuous by
Theorem~\ref{ctythm}. By Theorem \ref{attractorthm} and Lemma
\ref{intersectlemma},
\[
A=\lim_{L\rightarrow\infty}\overline{\bigcup_{j\geq L}\mathcal{F}^{j}(\left\{
x_{0}\right\}  )}.
\]
It follows that, for any $\varepsilon>0$, we can choose $K$ so that
\[
x_{k} \in\overline{\bigcup_{j\geq K}\mathcal{F}^{j}(\left\{
x_{0}\right\}  )} \subset A+\varepsilon
\]
for all $k\geq K$, as claimed. \vskip 2mm

We next show that, for any $\varepsilon>0$, there is an integer $K>0$
such that
\[
d_{\mathbb{H}}(A,\left\{  x_{k}\right\}  _{k=L}^{\infty} )<\varepsilon
\]
with probability one, for all $L\geq K$. This is equivalent to
the statement of the theorem.  To prove this, let 
 $\varepsilon> 0$. If $K$ is as specified in the paragraph above, then by (\ref{eq})
we have
$x_{L}\in A+\varepsilon$ for $L\geq K$. The attractor $A$, being compact,
is totally bounded. Let $\left\{  a_{q}:i=1,2,...,Q\right\}  $ be a set of
points such that $A\subset\cup_{q=1}^{Q} B(a_{i},\varepsilon/2)$, where
$B(a_{q},\varepsilon/2)$ is the ball of radius $\varepsilon/2$ centered at
$a_{q}$. Note that each $a_{q}$ and $Q$ depend on $\varepsilon$. By Lemma
\ref{uniformlemma} there is an integer $M$ such that, for each $x\in
A+\varepsilon$, there is $m < M$ such that $d_{\mathbb{H}}(A,\mathcal{F}%
^{m}(\left\{  x\right\}  )<\varepsilon/2.$ Hence
\[
d_{\mathbb{H}}(A,\mathcal{F}^{m}(\left\{  x_{L}\right\}  ))<\varepsilon/2
\]
for some integer $m < M$. Therefore there is a sequence of symbols \linebreak
$\sigma_{L+1}\sigma_{L+2}...\sigma_{L+m}$ such that $f_{\sigma_{L+m}}\circ
f_{\sigma_{L+m-2}}\circ...\circ f_{\sigma_{L+1}}(x_{L})\in B(a_{1}%
,\varepsilon/2)$. (We adopt the convention that the composition on the left
equals $x_{L}$ if $m=0$.) It follows that
\[
B(a_{1},\varepsilon/2)\; \bigcap\; \{x_{k}\}_{k=L}^{L+M-1} \neq\emptyset,
\]
or
\[
B(a_{1},\varepsilon/2)\subset\{x_{k}\}_{k=L}^{L+M-1} +\varepsilon.
\]
The probability that this event occurs, i.e., that the particular sequence
\linebreak$\sigma_{L+1}\sigma_{L+2}...\sigma_{L+m}$ is chosen, is greater than
$p^{M}$. By repeating this argument, we deduce that the probability that
\[
B(a_{q},\varepsilon/2)\subset\left\{  x_{k}\right\}  _{k=L+(q-1)M}^{L+qM-1}
+\varepsilon
\]
is greater than $p^{M} > 0$, for each $q\in\{1,2,...,Q\}$, regardless of
whether or not the preceding events occur. (That is not to say that the events
are independent.) It follows that the probability that all of these events
occur is greater than $p^{QM}$. Consider the event $E_{1}$ defined by
\[
\bigcup_{q=1}^{Q} B(a_{q},\varepsilon/2)\nsubseteq\{ x_{k}\}_{k=L}^{L+QM-1}
+\varepsilon
\]
The probability of $E_{1}$ is less than $(1-p^{QM})$. By a similar argument
the probability of the event $E_{r}, \, r\geq1$, defined by
\[
\bigcup_{q=1}^{Q}B(a_{q},\varepsilon/2)\nsubseteq\left\{  x_{k}\right\}
_{k=L+(r-1)QM} ^{L+rQM-1} \left\{  x_{k}\right\}  +\varepsilon
\]
is less than $(1-p^{QM})$, regardless of whether or not the previous events
\linebreak$E_{1},E_{2},...,E_{r-1}$ occurred, for $r=2,3, \dots$. It follows
that the probability of the event $E_{1}\cap E_{2}\cap\cdots\cap E_{r}$ is
less than $(1-p^{QM})^{r}$, for all $r=1,2,\dots$. This inequality holds
regardless of the fact that the $E_{r}$ are not independent. To simplify
notation, let $s = (1-p^{QM}) < 1$ so that
\[
pr (E_{1} \cap E_{2} \cap\cdots) \leq pr (E_{1} \cap E_{2} \cap\cdots\cap
E_{r}) = s^{r}%
\]
for all $r$. This implies that $pr(\cap_{r=1}^{\infty} E_{r}) = 0$. Hence,
with probability one, there is an $R$ such that
\[
\bigcup_{q=1}^{Q}B(a_{q},\varepsilon/2)\subset\left\{  x_{k}\right\}
_{k=L+(R-1)QM}^{L+RQM-1} +\varepsilon.
\]
Since $A\subset\cup_{q=1}^{Q}B(a_{q},\varepsilon/2)$ it follows that, with
probability one, there is $R$ such that
\[
A\subset\left\{  x_{k}\right\}  _{k=L+(R-1)QM}^{L+RQM-1} +\varepsilon .
\]
Because $L$ is an arbitrary integer greater than or equal to $K$,
\[
A \subset\left\{  x_{k}\right\}  _{L}^{\infty} +\varepsilon.
\]
for any $L\geq K$. But by (\ref{eq}) we also have
\[
\left\{  x_{k}\right\}  _{k=L}^{\infty} \subset A+\varepsilon.
\]
for any $L\geq K$. Hence, with probability one
\[
d_{\mathbb{H}} ( A, \left\{  x_{k}\right\}  _{k=L}^{\infty} )<\varepsilon.
\]
for any $L\geq K$.
\qed

It follows that
\[
\bigcap\limits_{K\geq1}\overline{\left\{
x_{k}\right\}_{k=K}^{\infty} }=\bigcap\limits_{K\geq1}\overline{\bigcup_{k\geq K}
\mathcal{F}^{k}(B)}
\]
almost surely,
for $x_{0}\in U$ and $B\in\mathbb{H}\left(  U\right) $, for example,
$B=\left\{  x_{0}\right\} $. We draw attention to this equality because it
seems suprising when $\mathcal{F}$ contains for than one function; the set
$\left\{  x_{k}\right\}_K^{\infty}$ seems sparse in comparison
to $\bigcup\limits_{k\geq K}\mathcal{F}^{k}(B)$.

\section{\label{examples}Examples}

\begin{example} The IFS $\mathcal{F}$ in this example has a unique 
attractor, yet each 
$f\in\mathcal{F}$ fails\label{ex1} to be a contraction with respect to
  any metric giving the same topology as the original metric.  With
  probability one, the chaos game applied to this example "draws a
  picture" of the attractor of the IFS.
\end{example}

If $\mathbb{X=\{(}x,y)\in\mathbb{R}^{2} :x^{2}+y^{2}=1\}$ and $d$ is
the Euclidean metric, then $(\mathbb{X},d)$ is a compact metric
space. Let $\mathcal{F}=(\mathbb{X};f_{1},f_{2})$ where
$f_{1}(x,y)=(x,y)$ and
$f_{2}(x,y)=(x\cos\alpha-y\sin\alpha, \, x\sin\alpha
+y\cos\alpha),$ where $\alpha/\pi$ is irrational. The
map $f_{1}$ is the identity map and $f_{2}$ is a rotation through
angle $\alpha$ anticlockwise about the origin. Since neither $f_{1}$
nor $f_{2}$ has a unique fixed point, it follows that there exists no
metric $\tilde{d}$ on $\mathbb{X}$ such that $(\mathbb{X},\tilde{d})$
is complete and either $f_{1}$ or $f_{2}$ is a contraction. On the
other hand, $\mathcal{F}$ has a unique attractor $A=\mathbb{X}$. To
see this, first note that $\mathcal{F}(\mathbb{X)=X}$.  Also, if
$(x_{0},y_{0})\in\mathbb{X}$ then
\[
\mathcal{F}^{k}(x_{0},y_{0})=  \{ f_2^j(x_{0},y_{0})\}_{j=0}^k
\]
for all $k$.
The right hand side is well known to converge in the Hausdorff metric, as $k$
tends to infinity, to $\mathbb{X}$. If follows that $\mathcal{F}^{k}(B)$
converges to $\mathbb{X}$ for all $B\in\mathbb{H(X)}$.  By
definition \ref{attractdef}, the IFS $\mathcal{F}$ has a unique attractor, 
namely
$\mathbb{X}$. By Theorem~\ref{chaosthm}, with probability one, the chaos game 
applied to this example "draws a picture" of the unit circle, the attractor of 
the IFS.

\begin{example} In the same spirit as Example\ref{ex1}, this is an IFS
on the real projective plane that possesses an attractor, but there is
no equivalent metric with respect to which the maps of the IFS are contractive.
Again, with probability one, the chaos game "draws a picture" of the
attractor of the IFS.
\end{example}

This example appears in \cite{BVW}.  The metric space is
$(\mathbb{R}P^{2},d)$, where $\mathbb{R}P^{2}$ denotes real projective
two-dimensional space, and $d$ denotes the round metric; see
\cite{BVW}.  Let $\mathcal{F}=\left(
\mathbb{R}P^{2};f_{1},f_{2}\right)$, where the pair of projective
transformations $f_{1}$ and $f_{2}$ are represented (acting on homogeneous
coordinates), respectively, by
the pair of matrices
\[%
\begin{pmatrix}
1 & 0 & 0\\
0 & 2 & 0\\
0 & 0 & 2
\end{pmatrix}
\qquad\text{and}\qquad%
\begin{pmatrix}
1 & 0 & 0\\
0 & 2\cos\alpha & -2\sin\alpha\\
0 & 2\sin\alpha & 2\cos\alpha
\end{pmatrix}
,
\]
and $\alpha/\pi$ is irrational. In terms of homogeneous coordinates
$(x,y,z)$, the attractor of $\mathcal{F}$ is the line $x=0$. This can be
proved by a similar argument to the one in Example \ref{ex1}.

Such non-contractive IFSs occur often in real projective IFS
theory. Theorem \ref{chaosthm} tells us that the chaos game algorithm
can always be applied to compute approximate pictures of attractors of
real projective IFSs.

\begin{example} This is a superfractal example.
\label{ex4} \end{example}

Let $(\mathbb{X},d_{\mathbb{X}})$ be a compact metric space.  Consider
the IFS \linebreak $\mathbb{F=(H(X)};\mathcal{F}_{1},\mathcal{F}_{2})$ where the
mappings $\mathcal{F}_{i}:\mathbb{H(X)\rightarrow H(X)}$ are the maps
corresponding to the IFSs $\mathcal{F}_{1}=(\mathbb{X};f_{11},f_{12})$
and $\mathcal{F}_{2} =(\mathbb{X};f_{21},f_{22})$ - as explained
directly after definition~\ref{defIFS}.  The space $\mathbb{H(X)}$ is
endowed with the Hausdorff metric $d_{\mathbb{H(X)}}$ derived from the
Euclidean metric $d_{\mathbb{X}}$ on $\mathbb{X}$. In particular,
since the functions $f_{ij}$ are continuous, the mappings
$\mathcal{F}_{i}:\mathbb{H(X)\rightarrow H(X)}$ are indeed continuous
and $\mathbb{F}$ is a well-defined IFS. By Theorem \ref{ctythm},
$(\mathbb{H(X)},d_{\mathbb{H(X)}})$ is a compact, and hence normal,
metric space. Again invoking Theorem \ref{ctythm}, it follows that
$(\mathbb{H(H(X))},d_{\mathbb{H(H(X))}})$ is also a compact, and hence
normal metric space. It follows that $\mathbb{F}:$
$\mathbb{H(H(X))\rightarrow}$ $\mathbb{H(H(X))}$ defined by%
\[
\mathbb{F(B)=}\mathcal{F}_{1}(\mathbb{B)\cup}\mathcal{F}_{2}(\mathbb{B)}
\]
for all $\mathbb{B\in H(H(X))}$.
is continuous with respect to the Hausdorff-Hausdorff metric
$d_{\mathbb{H(H(X))}}$, where
\[
\mathcal{F}_{i}\mathbb{(B)=\{}\mathcal{F}_{i}(B):B\in\mathbb{B\}\subset
H(H(X))}
\]
for all $\mathbb{B}\in\mathbb{H(H(X))}$.
Consistent with the terminology in \cite{barnsuper1},
\cite{barnsuper2}, \cite{barnsuper3}, an attractor of such an
IFS $\mathbb{F}$ is called a {\it superfractal}. A superfractal thus
consists of a set $\mathbb{A}$ whose elements are compact
subsets of $\mathbb{X}$. The set $\mathbb{A}$ is itself compact in the
Hausdorff-Hausdorff metric. Theorem 3 tells us that such an attractor
can be computed by means of the chaos game, regardless of whether or
not there exists a metric such that the constituent mappings are
contractive. Many examples of approximations to sets 
belonging to superfractals, calculated using the chaos game algorithm,
are illustrated in \cite{barnsuper1}. These calculations were made,
without any special care over whether or not the random selections
were precisely i.i.d.  or even whether or not the constituent IFS's
were contractive. This paper justifies such apparent carelessness.


\begin{thebibliography}{99}                                                                                               %
\bibitem {BVW}M. F. Barnsley, A. Vince, D. C. Wilson, Real projective iterated
function systems, \textit{Preprint (2010).}

\bibitem {barnsley}M. F. Barnsley, \textit{Fractals Everywhere}, Academic
Press, Boston, MA, 1988.

\bibitem {barnsleydemko}M. F. Barnsley and S. G. Demko,\textit{\ }Iterated
function systems and the global construction of fractals, \textit{Proc. Roy.
Soc. London Ser. A} \textbf{399} (1985) 243--275.

\bibitem {barnplacedep}M. F. Barnsley, S. G. Demko, J. H. Elton, J. S.
Geronimo, Invariant measures arising from iterated function systems with
place-dependent probabilities, \textit{Ann. Inst. H. Poincar\'{e}}, \textbf{24
}(1988) 367-394.

\bibitem {barnsuper1}Michael F. Barnsley, \textit{Superfractals}, Cambridge
University Press, Cambridge, 2006.

\bibitem {barnsuper2}M. F. Barnsley, J. Hutchinson, \"{O}. Stenflo, V-variable
fractals: fractals with partial self similarity, \textit{Advances in
Mathematics}, \textbf{218 (}2008) 2051-2088.

\bibitem {barnsuper3}M. F. Barnsley, J. Hutchinson, \"{O}. Stenflo,
$V$-variable fractals, dimension results, \textit{Forum Mathematicum, }to appear.

\bibitem {berger}Marc A. Berger, Random affine iterated function systems:
curve generation and wavelets, \textit{SIAM Review, }\textbf{34} (1992) 361-385.

\bibitem {mcgehee}R. P. McGehee, T. Wiandt, Conley decomposition for closed
relations\textit{, Difference Equations and Applications, \textbf{12 }(2006)
1-47.}

\bibitem {dudley}R. M. Dudley, \textit{Real Analysis and Probability,
}Wadsworth \&\ Brooks/Cole, Pacific Grove, California, 1989.

\bibitem {JEergodic}J. H. Elton, An ergodic theorem for iterated
maps\textit{,} \textit{Ergodic Theory Dynam. Systems,} \textbf{7} (1987) 481--488.

\bibitem {edgar}Gerald A. Edgar, \textit{Measure, Topology, and Fractal
Geometry}, Springer-Verlag, New York, 1990.

\bibitem {jaroszewska}Joanna Jaroszewska, Iterated function systems with
continuous place dependent probabilities,\textit{ Universitatis Iagellonicae
Acta Mathematica, }\textbf{40 }(2002) 137-146.

\bibitem {henrikson}Jeff Henrikson, Completeness and total boundednesss of the
Hausdorff metric,\textit{ MIT Undergraduate Journal of Mathematics, \textbf{1
}}(1999)\textit{ 69-79.}

\bibitem {vrscay1}B. Forte and\ F. Mendivil, A classical ergodic property for
IFS: A simple proof, \textit{Ergodic Theory Dynam. Systems} \textbf{18} (1998) 609--611.

\bibitem {hutchinson}J. E. Hutchinson, Fractals and self-similarity\textit{,
Indiana Univ. Math. J.} \textbf{30 }(1981) 713--747.

\bibitem {jamison}Benton Jamison, Asymptotic stability of successive iterates
of continuous functions under a Markov operator, \textit{Journal of
mathematical analysis and applications, }\textbf{9 }(1964) 203-214.

\bibitem {lesniak}Krzysztof Le\'{s}niak, Stability and invariance of
multivalued iterated function systems, \textit{Math. Slovaca, }\textbf{53
}(2003) 393-405.

\bibitem {onicescu}O. Onicescu and G. Mihok, Sur les cha\^{\i}nes de
variables\ statistiques, \textit{Bull. Sci. Math. de France} \textbf{59}
(1935) 174-192.

\bibitem {stenflo3}\"{O} Stenflo, Uniqueness of invariant measures for
place-dependent random iterations of functions\textit{, }in\textit{\ Fractals
in Multimedia, }M. F. Barnsley, D. Saupe, and E. R. Vrscay, eds., IMA Volumes
in Mathematics and its Applications, vol. 132, Springer-Verlag, New York, 2002.

\bibitem {werner}Ivan Werner, Ergodic theorem for contractive Markov systems,
\textit{Nonlinearity, }\textbf{17} (2004) 2303-2313.

\bibitem {S}Robert Scealy, \textit{$V$-variable fractals and interpolation},
Ph.D. Thesis, Australian National University, 2008.
\end{thebibliography}
\end{document}